\documentclass[graybox]{svmult}

\usepackage[backend=biber,
            style=numeric-comp, % numerico + comprime intervalli [3–5]
            sorting=none,       % ordine di citazione
            doi=true, url=false]{biblatex}
\usepackage{csquotes}        % good practice with biblatex

\usepackage{hyperref}
\usepackage{type1cm}        
\usepackage{makeidx}         
\usepackage{graphicx}        
\usepackage{multicol}        
\usepackage[bottom]{footmisc}\usepackage{caption}
\usepackage{algorithm,algpseudocode}
\usepackage{pgfplotstable}
\usepackage{pgfplots}
\usepackage{array,multirow,graphicx}
\usepackage{pgfplots}
\pgfplotsset{compat=newest}
    \usepgfplotslibrary{groupplots}
    \usetikzlibrary{matrix}
\DeclareMathAlphabet{\pazocal}{OMS}{zplm}{m}{n}    
\pgfplotsset{compat = 1.3}
\usepackage{wrapfig}

\usepackage{xcolor}
\definecolor{TR}{HTML}{6a4c93}
\definecolor{SN}{HTML}{c9182c}

\usepackage{mathtools}
\mathtoolsset{showonlyrefs}

\usepackage{booktabs}

\definecolor{TLTR0}{rgb}{0.63,0.12156,0.48236} 
\definecolor{TLTR1}{HTML}{ff924c}
\definecolor{TLTR2}{HTML}{ffd700}
\definecolor{TLTR3}{HTML}{8ac926}
\definecolor{TLTR4}{HTML}{1982c4}    

\newcommand{\rH}{\mathrm{H}}
\newcommand{\rL}{\mathrm{L}}
\newcommand{\trsp}{\top}

\def \0{\mathbf{0}}			% 

\def \R{{\mathbb{R}}}			% real num
\def \N {\mathbb{N}}			% numbers 
\def \I{{\mathbf{I}}}			% subdivision matrix
\def \0{\vec{0}}			% 

\def \R{\mathbb{R}}			% real numbers
\def \N {\mathbb{N}}			% nat. numbers

\def \I{\mathbf{I}}			% subdivision matrix

\usepackage{newtxtext}       
\usepackage[varvw]{newtxmath}       

\algnewcommand\algorithmiconput{\textbf{Constants:}}
\algnewcommand\algorithmicinput{\textbf{Input:}}
\algnewcommand\algorithmicoutput{\textbf{Output:}}
\algnewcommand{\algorithmicgoto}{\textbf{go to}}%

\algnewcommand\Constants{\item[\algorithmiconput]}
\algnewcommand\Input{\item[\algorithmicinput]}%
\algnewcommand\Output{\item[\algorithmicoutput]}%
\algnewcommand{\Goto}[1]{\algorithmicgoto~\ref{#1}}%

\makeindex

\begin{document}
%Magical trust region methods with...
\title*{Trust-Region Methods with Low-Fidelity Objective Models}%A Multifidelity Trust-Region Method for Classifier Training %A Trust-Region Method with Feature-Based Low-Fidelity Directions for Classifier Training
\author{Andrea Angino\orcidID{0009-0000-8525-375X},\\
Matteo Aurina\orcidID{0009-0002-2654-7488},\\
Alena Kopaničáková\orcidID{0000-0001-8388-5518},\\
Matthias Voigt\orcidID{0000-0001-8491-1861},\\
Marco Donatelli\orcidID{0000-0001-7958-9126}, and\\
Rolf Krause\orcidID{0000-0001-5408-5271}.}
\authorrunning{Andrea Angino et al.}
\institute{
Andrea Angino \at UniDistance Suisse; \email{andrea.angino@unidistance.ch}
\and Matteo Aurina \at Università dell’Insubria; \email{Matteoaurina@gmail.com}
\and Alena Kopaničáková \at Toulouse INP-ENSEEIHT, IRIT, ANITI; \email{alena.kopanicakova@toulouse-inp.fr}
\and Matthias Voigt \at UniDistance Suisse; \email{matthias.voigt@fernuni.ch}
\and Marco Donatelli \at Università dell’Insubria; \email{marco.donatelli@uninsubria.it}
\and Rolf Krause \at King Abdullah University of Science and Technology (KAUST), UniDistance Suisse; \email{rolf.krause@kaust.edu.sa}
}

\maketitle
\abstract
{
We introduce two multifidelity trust-region methods based on the Magical Trust Region (MTR) framework. MTR augments the classical trust-region step with a secondary, informative direction. In our approaches, the secondary ``magical'' directions are determined by solving coarse trust-region subproblems based on low-fidelity objective models.
The first proposed method, Sketched Trust-Region (STR), constructs this secondary direction using a sketched matrix to reduce the dimensionality of the trust-region subproblem. 
The second method, SVD Trust-Region (SVDTR), defines the magical direction via a truncated singular value decomposition of the dataset, capturing the leading directions of variability. Several numerical examples illustrate the potential gain in efficiency.

}

\section{Introduction}
We consider large-scale unconstrained optimization problems arising from data-driven applications, particularly those encountered in supervised machine learning. Specifically, we focus on binary classification tasks, where the training dataset is given by
\[
\mathcal{D} = \{(x_i, y_i) \in \mathbb{R}^n \times \mathcal{Y} \mid i = 1, \dots, q\},
\]
with \(x_i \in \mathbb{R}^n\) representing the $i$-th feature vector and \(y_i \in \mathcal{Y}= \{-1, 1\}\) the corresponding labels. Furthermore, we collect the feature vectors into the data matrix
\[
X = [x_1, x_2, \dots, x_q] \in \mathbb{R}^{n \times q},
\]
which will be used to define low-fidelity directions in the methods presented below.

%The learning task is formulated as the minimization of an empirical risk functional, where the objective is a finite sum of individual loss terms evaluated over the dataset, i.e.,

The learning task is formulated as an empirical risk  minimization, where the objective is a finite sum of individual loss terms evaluated over the dataset, i.e.,
\begin{equation}\label{eq:Problem}
    \min_{w \in \mathbb{R}^n} f(w) = \frac{1}{q} \sum_{i=1}^q \ell\big(w; x_i, y_i\big),
\end{equation}
where the optimization variable \(w \in \mathbb{R}^n\) denotes classifier weights and \(\ell(w; x_i, y_i)\) is a smooth loss function (e.g., logistic, squared hinge, cross-entropy loss) measuring the misfit between the model prediction and the label for the $i$-th data point. We assume that \(f\) is bounded from below and twice continuously differentiable with respect to \(w\in\mathbb{R}^n\). Both the dimensionality of the parameter vector \(n\) and the dataset size \(q\) can be large in modern applications. 

While first-order methods, such as stochastic gradient descent (SGD)~\cite{RobbinsMonro1951} and adaptive schemes like Adam~\cite{KingmaBa2015}, are widely used in large-scale machine learning due to their simplicity and low per-iteration computational cost, they often suffer from slow convergence and sensitivity to hyperparameter tuning, particularly in nonconvex settings. Among alternative methods, trust-region (TR) algorithms~\cite{ConnTRmethods} construct a local model of the objective function at each iteration and solve a constrained subproblem to determine the search direction, guaranteeing global convergence to a first-order critical point~\cite{curtis2017trust}. Classical enhancements to TR include two-step variants tailored to structured problems~\cite{ConnMTR}.% and they interface naturally with domain decomposition ideas.% and second-order correction mechanisms to improve step quality~\cite{Fletcher1982}.

In this work, building on ideas from multifidelity optimization~\cite{Forrester2007,PeherstorferSurvey2018}, we introduce two multifidelity trust-region methods inspired by the \emph{magical} trust-region (MTR) framework~\cite[Section 10.4.1]{ConnTRmethods}, which augment classical TR steps with an additional ``magical" direction aimed at accelerating convergence, see the recent two-level TR method in the same framework~\cite{AnginoKopanicakovaKrause2024TLTR}. Traditionally, the MTR framework assumes the availability of an oracle that provides enhanced directions, improving upon those obtained from the standard model.

Our first method, called Sketched Trust-Region (STR), constructs the secondary direction by
sketching the data matrix \(X\) at every iteration, thereby reducing the dimensionality of the trust–region subproblem. In contrast to classical sketched optimization methods that rely entirely on the reduced model~\cite{newtonSketch, sketchingGD, Cartis}, STR employs the sketch matrix to only generate a corrective low-fidelity direction that enhances the full-space TR step.

The second method, named SVD Trust-Region (SVDTR), defines the magical direction via a truncated SVD of the data matrix \(X\), retaining the leading \(t\) singular vectors to form the feature projector. This captures the dominant directions of variability in the dataset, which is particularly effective when the singular values decay rapidly. 

Viewed through the lens of domain decomposition, STR provides an algebraic coarse correction via
on-the-fly feature aggregation, whereas SVDTR supplies a spectral coarse space from the
dominant singular vectors of \(X\).

%In both methods, the low-fidelity direction influences convergence speed but not convergence itself, allowing for dimensionality reduction. 
Our goal is to apply these approaches to classification tasks in machine learning, where the balance between cost and accuracy is critical. By incorporating data-driven low-fidelity models into each TR step, we aim to improve the efficiency of training procedures.

%This manuscript is organized as follows.  
%In Section~\ref{Sec:DMTR}, we present the STR and SVDTR algorithms in detail.  
%In Section~\ref{Sec:NumericalEX}, we evaluate their performance on neural network training problems and compare them with standard methods.

\section{Magical TR with low-fidelity directions}
\label{Sec:DMTR}
%In this section, we present two multifidelity trust-region methods, STR and SVDTR, designed for large-scale supervised learning tasks. Both methods introduce a composite search direction that combines information from the high-fidelity objective with a data-driven low-fidelity model. The primary direction is computed by solving a standard trust-region subproblem on the full objective~\eqref{eq:Problem}. The secondary direction is obtained from a reduced objective constructed via feature extraction on the training data: STR uses randomized sketching to define a low-dimensional subspace, while SVDTR employs a truncated singular value decomposition of the dataset to capture the leading directions of variability.

%\subsection{The STR and SVDTR algorithms}
Following the MTR framework~\cite{ConnTRmethods}, both STR and SVDTR are initialized with an initial guess \( w_0 \in \mathbb{R}^n \). At the \(k\)-th iteration, the algorithms first compute a high-fidelity step \( p_k^\rH \) by approximately solving the trust-region subproblem
\begin{align}
\begin{split}
& \min_{p_k^\rH \in \mathbb{R}^n} m_k^\rH\big(p_k^\rH\big) := f(w_k) + \left\langle \nabla f(w_k), p_k^\rH \right\rangle + \frac{1}{2} \left\langle p_k^\rH, \nabla^2 f(w_k) p_k^\rH \right\rangle, \\
& \text{subject to} \quad \left\|p_k^\rH\right\| \leq \Delta_k,
\end{split}
\label{eq:modelfull}
\end{align}
where \( m_k^\rH \) is the quadratic model of the full objective around \( w_k \) and \( \Delta_k > 0 \) is the trust-region radius controlling the step size.

The secondary, low-fidelity direction is then computed around the intermediate iterate \( w_{k+1/2} := w_k + p_k^\rH \). Both methods define a low-dimensional objective
\begin{equation}
f_k^\rL(\tilde w) := \frac{1}{q} \sum_{i=1}^q \tilde \ell(\tilde w; \tilde x_i, y_i),
\label{eq:sketched_obj}
\end{equation}
where \( \tilde x_i = S_k x_i \) are the reduced feature vectors obtained via a projection matrix \( S_k \in \mathbb{R}^{t \times n} \) with \( t \ll n \). 
\begin{itemize}
    \item In STR, \(S_k\) is a randomized sketching matrix (e.g., Gaussian) used to compress the dataset while approximately preserving its geometric structure~\cite{Lindenstrauss}.
  \item In SVDTR, \(S_k\) is a fixed matrix across iterations (the subscript \(k\) is retained for consistency in the algorithm notation), can be pre-computed, and is constructed from the leading left singular vectors of \(X\).
\end{itemize}

A second-order model \( m_k^\rL \) of \( f_k^\rL \) is built around \( \tilde w_{k+1/2} = S_k w_{k+1/2} \), and the low-fidelity step \( p_k^\rL \in \mathbb{R}^t \) is computed by solving the trust-region subproblem
\begin{align}
\begin{split}
& \min_{p_k^\rL \in \mathbb{R}^t} m_k^\rL\big(p_k^\rL\big) := f_k^\rL(\tilde w_{k+1/2}) + \left\langle \nabla f_k^\rL(\tilde w_{k+1/2}), p_k^\rL \right\rangle + \frac{1}{2} \left\langle p_k^\rL, \nabla^2 f_k^\rL(\tilde w_{k+1/2}) p_k^\rL \right\rangle, \\
& \text{subject to} \quad \left\|p_k^\rL\right\| \leq \Delta_k.
\label{eq:modellow}
\end{split}
\end{align}

The reduced step is then lifted to the full space via \( S_k^\trsp p_k^\rL \) and incorporated into the update only if it decreases the original objective, i.e.,
\begin{equation}
f\big(w_k + p_k^\rH + \alpha_k S_k^\trsp p_k^\rL\big) < f\big(w_k + p_k^\rH\big),
\label{eq:conditionstep}
\end{equation}
where \( \alpha_k > 0 \) may be fixed or chosen via a line search strategy; otherwise we set \( p_k^\rL = 0 \). When \( p_k^\rL = 0 \), the algorithm reduces to a classical trust-region method, ensuring global convergence. 

The effectiveness of the composite step \( p_k := p_k^\rH + \alpha_k S_k^\trsp p_k^\rL \) is measured by a trust-region ratio
\begin{equation}
\varrho_k := \frac{f(w_k) - f(w_k + p_k)}{m_k^\rH(w_k) - m_k^\rH\big(w_k + p_k^\rH\big) + f\big(w_k + p_k^\rH\big) - f(w_k + p_k)},
\label{eq:rho_dmtr}
\end{equation}
which determines step acceptance and trust-region radius updates. Thus, the low-fidelity step may improve acceptance of steps that would be rejected by standard TR, accelerating convergence. The complete procedure for both methods is summarized in Algorithm~\ref{ALG:DMTR}, where the only difference lies in the construction of the low-fidelity feature projection.

\begin{algorithm}[t]
\footnotesize
\caption{\footnotesize Sketched Trust-Region Method}
\label{ALG:DMTR}
\begin{algorithmic}[1]
\Require {$ f: \R^n \rightarrow \R,  w_0 \in \R^n,  \Delta_{0} \in \R^+, t < n \in \N$}
\Ensure Minimizer $w^*$ of $f$
\Constants { $0 < \eta_1 \leq \eta_2 < 1, \  0 < \gamma_1 \leq  \gamma_2 < 1$}
\State $k:=0$
\While{not converged}
\State $p_k^\rH := \underset{\left\| p \right\| \leq \Delta_k}{\text{argmin}}  \;m_k^\rH\left(p\right)$ \Comment{Obtain full-space search direction}
\State $w_{k+1/2} :=w_k + p_k^\rH$ 
\State Construct $S_k$ via sketching  \Comment{For SVDTR: $S_k$ is precomputed}
\State $\tilde X :=S_k X$
\State $ p_k^\rL :=  \underset{\| {\tilde p} \| \leq \Delta_k}{\text{argmin}} \;m_k^\rL\left({\tilde p}\right) $ \Comment{Obtain subspace search-direction }
\State 
\vspace{-0.4cm}
\Comment{Assess the quality of the subspace step}
\begin{align*}
& p_{k}^\rL := 
\begin{cases}
{p}_k^\rL, \quad &\text{if} \ f\left(w_k+p_k^\rH + \alpha_k S_k^\trsp {p}_k^\rL\right) < f\left(w_k + p_k^\rH \right)    \\
0, \quad & \text{otherwise} \hspace{6.75cm} \textcolor{white}{.}
\end{cases}  &
\end{align*} 
\State Evaluate $\varrho_k$ as in~\eqref{eq:rho_dmtr} \Comment{Assess the quality of the composite trial step}
\vspace{-0.3cm}
\begin{align*}
 & w_{k+1} := 
\begin{cases}
w_{k} + p_k^\rH + \alpha_k S_k^\trsp p_{k}^\rL, \quad &\text{if} \ \varrho_k > \eta_1,   \\
w_{k}, \quad & \text{otherwise,}
\end{cases} &
& \Delta_{k+1} := 
\begin{cases}
[\Delta_k, \infty), \quad &\text{if} \ \ \varrho_k \geq \eta_2,   \\
[\gamma_2 \Delta_k, \Delta_k], \quad &\text{if} \ \ \varrho_k \in [\eta_1, \eta_2),   \\
[\gamma_1 \Delta_k, \gamma_2 \Delta_k], \quad &\text{if} \ \ \varrho_k < \eta_1  
\end{cases} 
\end{align*}
\State $k := k+1$
\EndWhile
\State \Return $w^* := w_{k}$
\end{algorithmic}
\end{algorithm}
%\subsection{The computational cost of STR and SVDTR}  
In both methods, the high-fidelity step \(p_k^\rH\) is computed approximately (e.g., via a few Steihaug-Toint CG iterations or using the Cauchy point).  The additional computational effort compared to classical trust-region methods arises from two main tasks: constructing the low-fidelity model and solving the corresponding reduced trust-region subproblem.

\section{Numerical examples}
\label{Sec:NumericalEX}
We evaluate the proposed algorithms, STR and SVDTR, on binary classification problems for the empirical-risk formulation \eqref{eq:Problem}, and compare them against the classical full-space TR baseline. We consider two objective functions, namely
\begin{equation}\label{eq:losses}
\begin{aligned}
f_{\mathrm{LL}}(w) &= \frac{1}{q}\sum_{i=1}^q \log \left(1+\mathrm{e}^{-y_i\langle w,x_i\rangle}\right)
                  + \tfrac{\lambda}{2}\|w\|_2^2,\\%[4pt]
f_{\mathrm{LS}}(w) &= \frac{1}{q}\sum_{i=1}^q \Biggl(y_i - \frac{\mathrm{e}^{\langle w,x_i\rangle}}{1+\mathrm{e}^{\langle w,x_i\rangle}}\Biggr)^{\!2}
                  + \tfrac{\lambda}{2}\|w\|_2^2,
\end{aligned}
\end{equation}
with regularization parameter \(\lambda = 1/q\).

Unless stated otherwise, each run is terminated when the Euclidean norm of the full gradient satisfies \({\|\nabla f(w)\|}_2 \le 10^{-6}\) or when a predefined iteration/wall-clock time budget is reached. %All plots report the full-gradient norm \(\|\nabla f(w_k)\|_2\) versus iterations; wall-clock time is shown when informative.
The datasets are drawn from the LIBSVM repository (\texttt{Australian} (621 samples, 14 features), \texttt{Mushroom} (6,499 samples, 112 features), \texttt{Gisette} (6,000 samples, 5,000 features)).\footnote{\url{https://www.csie.ntu.edu.tw/~cjlin/libsvmtools/datasets/binary.html}} %Table~\ref{tab:dataset_characteristics} lists the number of samples \(q\) and features \(n\).
The methods are implemented in \texttt{Python} using \texttt{PyTorch}~2.8.0~\cite{paszke2019pytorch}. 
%Our code runs on both CPU- or GPU-based architectures (when available); 
All reported results were obtained on Windows~64-bit with an AMD Ryzen~7~5700G CPU (3.80\,GHz) and 16\,GB RAM (CPU-only).

The high-fidelity TR subproblems in \eqref{eq:modelfull} are approximately solved with the Steihaug-Toint conjugate gradient (ST-CG) method with two inner iterations for the \texttt{Australian} and \texttt{Mushroom} datasets. For the \texttt{Gisette} dataset, we either use ST-CG with $25$ inner iterations or the Cauchy-point (CP) solver.

% \noindent
% \begin{minipage}[t]{0.49\textwidth}
%   \centering
%   \captionof{table}{Characteristics of the datasets.}
%   \label{tab:dataset_characteristics}
%   \begin{tabular}{ccc}
%     \toprule
%     dataset & $q$ (samples) & $n$ (features) \\
%     \midrule
%     Australian & 621   & 14 \\
%     Mushroom   & 6{,}499 & 112 \\
%     Gisette    & 6{,}000 & 5{,}000 \\
%     \bottomrule
%   \end{tabular}
% \end{minipage}\hfill
% \begin{minipage}[t]{0.49\textwidth}
%   For STR and SVDTR, the low-fidelity subproblems in \eqref{eq:modellow},
%   posed in a reduced space of dimension \(t\), are solved by ST-CG with at most \(t\)
%   inner iterations, so the reduced directions are effectively accurate once lifted to
%   the full space. In STR, the sketch matrix \(S\in\mathbb{R}^{t\times n}\) has i.i.d. entries
%   \(\mathcal{N}(0,t^{-1})\). In SVDTR, the reduced
% \end{minipage}
 For STR and SVDTR, the low-fidelity subproblems in \eqref{eq:modellow},
  posed in a reduced space of dimension \(t\), are solved by ST-CG with at most \(t\)
  inner iterations, so the reduced directions are effectively accurate once lifted to
  the full space. In STR, the sketch matrix \(S\in\mathbb{R}^{t\times n}\) has i.i.d. entries drawn from
  \(\mathcal{N}(0,t^{-1})\). In SVDTR, the reduced space is the span of the top \(t\) left singular vectors of \(X\).

\def\resbase{results}

\begin{figure}[bt]
    \centering
    \includegraphics[width=1\linewidth]{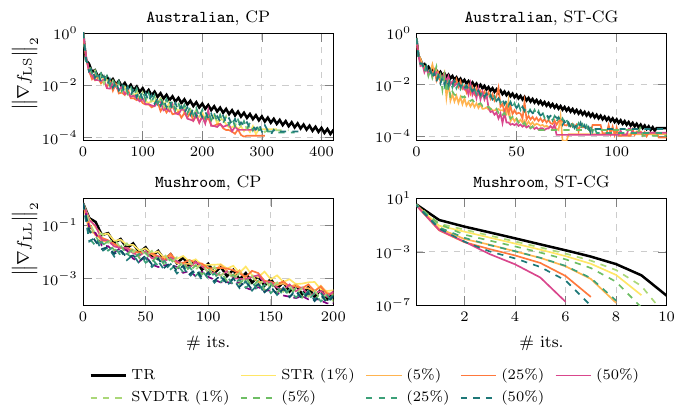}
    \caption{\footnotesize
Convergence histories of TR (solid black), STR (solid), and SVDTR (dashed) for solving~\eqref{eq:Problem}.  Top: \texttt{Australian} with $f_{\mathrm{LS}}$ using CP (left) and ST–CG (right). Bottom: \texttt{Mushroom} with $f_{\mathrm{LL}}$ under the same full–space solvers. Legend entries for STR/SVDTR indicate the reduced dimension $t$ as a percentage of the feature dimension $n$.}
    \label{fig:Austr_Mush}
\end{figure}

Figure~\ref{fig:Austr_Mush} displays the evolution of ${\|\nabla f(w_k)\|}_2$ as a function of outer iterations.
The top/bottom row reports the results for the \texttt{Australian}/\texttt{Mushroom} dataset with the $f_{\mathrm{LS}}$/$f_{\mathrm{LL}}$ loss using CP (left) and ST–CG (right). 
% The top row reports the results for the \texttt{Australian} dataset with the least–squares loss $f_{\mathrm{LS}}$ using CP (left) and ST–CG (right). The bottom row reports the results for the \texttt{Mushroom} dataset with the logistic loss $f_{\mathrm{LL}}$ under the same full–space solvers.
Wall–clock times are not reported, as for these small–scale datasets, all solvers complete within negligible runtime, making iteration counts the most meaningful comparison.
Across all configurations, augmenting the full–space step with a reduced–space direction yields a systematic reduction in the number of outer iterations required to attain comparable gradient norms.
The improvement is monotone with $t$ and is most pronounced when the full–space subproblems are solved by ST–CG; CP exhibits the same qualitative trend, albeit with smaller margins.
These observations substantiate the effectiveness of the proposed two–direction TR framework in accelerating convergence.

\begin{figure}[bt]
    \centering
     \includegraphics[width=0.95\linewidth]{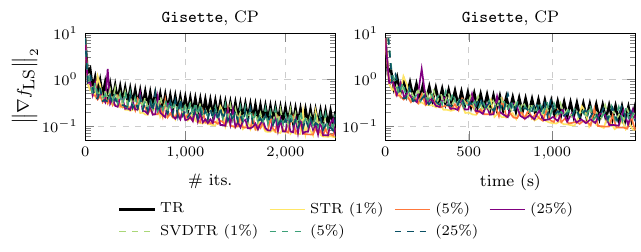}
    \caption{\footnotesize Convergence histories of TR (solid black), STR (solid), and SVDTR (dashed) for solving~\eqref{eq:Problem} with \(f_{\mathrm{LS}}\), all using CP. Left: \([\lVert \nabla f_{\mathrm{LL}} \rVert]_{2}\) versus iteration count; right: \({\lVert \nabla f_{\mathrm{LS}} \rVert}_{2}\) versus wall-clock time (s). Legend entries for STR/SVDTR indicate the reduced dimension \(t\) as a percentage of  \(n\).
}
    \label{fig:gisetteLS}
\end{figure}

\begin{figure}[bt]
    \centering
    \includegraphics[width=1\linewidth]{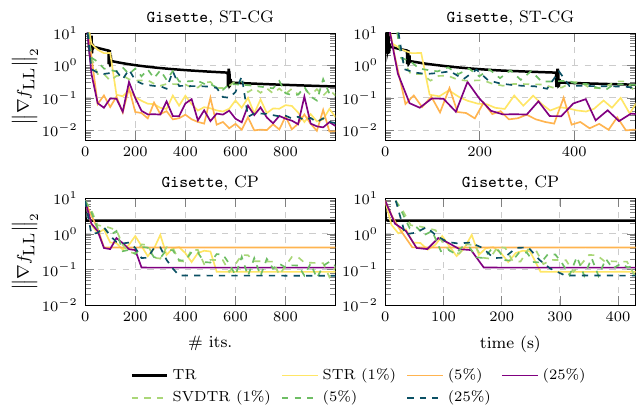}
     \caption{\footnotesize Convergence histories on the \texttt{Gisette} dataset of TR (solid black), STR (solid), and SVDTR (dashed) for solving~\eqref{eq:Problem} with $f_{\mathrm{LL}}$. Top: ST–CG full–space solver. Bottom: CP full–space solver.  Legend entries for STR/SVDTR indicate the reduced dimension $t$ as a percentage of $n$.
     }
    \label{fig:gisette_plots}
\end{figure}

We proceed by testing our methods on the high-dimensional \texttt{Gisette} dataset.
Figures~\ref{fig:gisetteLS} and
~\ref{fig:gisette_plots}  report the decay of the full-gradient norm \({\|\nabla f(w_k)\|}_2\) versus outer iterations and wall–clock time in seconds  for TR, STR, and SVDTR at multiple reduced dimensions \(t\).
In all cases, augmenting the full–space step with a reduced–space direction markedly lowers the iteration counts relative to TR, with a monotone trend as \(t\) increases. Overall, {SVDTR} tends to excel for \(f_{\mathrm{LL}}\) once the reduced dimension \(t\) is sufficiently large for the subspace to capture the dataset structure, whereas {STR} provides robust preprocessing-free improvements.

%\section{Summary}
%On small/medium datasets, {STR} and {SVDTR} reduce outer iterations with negligible time reduction. On larger datasets, both methods can improve both the  iteration count and wall-clock time for suitable \(t\), demonstrating the benefit of incorporating a low-fidelity direction within TR.

%\vspace{-0.2cm}
\section*{Code and data availability}
The code and data used to produce the numerical results is under available at
\begin{center}
  \url{https://doi.org/10.5281/zenodo.17473878}.
\end{center}
\begin{acknowledgement}
The work of A.K. benefited from the AI Interdisciplinary Institute ANITI, funded by the  France 2030 program under Grant Agreement No. ANR-23-IACL-0002. The research of A.A. and M.V. was funded in part by the Swiss National Science Foundation (SNSF) grant No. 224943.
\end{acknowledgement}

%\vspace{-1.2cm}

\printbibliography

\end{document}